\documentclass[preprint,12pt]{amsart}

\usepackage{amssymb}
\usepackage{amsthm}
\usepackage{amsmath}
\usepackage{mathrsfs}

\newtheorem*{theorem*}{Theorem}
\newtheorem{lemma}{Lemma} 

\theoremstyle{definition} 
\newtheorem{definition}{Definition}
\newtheorem*{problem*}{Problem}

\newcommand{\+}{\mathop{+_{\scriptscriptstyle 2}}}

\begin{document}

{
\renewcommand*{\thefootnote}{$\star$}

\title{A Non-$\mathbb R$-Factorizable 
Product\\ of $\mathbb R$-Factorizable Groups}
\footnotetext[0]{This work was financially supported by the Russian Science Foundation, grant~22-11-00075-P.} 

}

\author{Ol'ga Sipacheva}

\email{osipa@gmail.com}

\address{Department of General Topology and Geometry, Faculty of Mechanics and  Mathematics, 
M.~V.~Lomonosov Moscow State University, Leninskie Gory 1, Moscow, 199991 Russia}

\begin{abstract}
An example of two zero-dimensional $\mathbb R$-factorizable groups whose product is not $\mathbb R$-factorizable 
is constructed. One of these groups is second countable  and the other Lindel\"of to any finite power. 
\end{abstract}



\maketitle

\begin{definition}[\cite{Tk1}]
A topological group $G$ is said to be \emph{$\mathbb R$-factorizable} if, for every continuous function 
$f\colon G\to \mathbb R$, there exists a continuous homomorphism $h\colon G \to H$ to a second countable 
topological group $H$ and a continuous function $g\colon H\to \mathbb R$ such that $f = g \circ h$. 
\end{definition}

The study of $\mathbb R$-factorizable groups goes back to the work of Pontryagin, who proved the 
$\mathbb R$-factorizability of compact groups \cite[Example~37]{Pontryagin} (see also 
\cite[Theorem~8.1.1]{AT}), although the notion was explicitly introduced only as late as 1991 by Tkachenko in 
 \cite{Tk1}. In the same paper Tkachenko asked the question of whether or not the $\mathbb R$-factorizability of 
groups is preserved by finite products \cite[Problem~4.1]{Tk1}; versions of this question (some of which 
still remain unanswered) can be found in~\cite{AT}. 

The first examples of $\mathbb R$-factorizable groups $G$ and $H$ for which $G\times H$ is not 
$\mathbb R$-factorizable were given by this author~\cite{preprint} and, independently, Reznichenko~\cite{EAR}. 
All these examples were Lindel\"of and had some additional properties (for example, Reznichenko constructed 
a pair of Lindel\"of groups whose product was not pseudo-$\aleph_1$-compact and another pair of 
Lindel\"of groups  whose product was separable and contained a closed discrete subspace of cardinality 
$2^\omega$). In this paper, we construct two zero-dimensional 
$\mathbb R$-factorizable groups $G_1$ and $G_2$ such that $G_2$ is second countable, $G_1^n$ is Lindel\"of 
for any positive integer $n$, and $G_1\times G_2$ is not $\mathbb R$-factorizable, thereby solving  
Problems~8.5.2, 8.5.4, and 8.5.6 and one half of Problem~8.5.5 in \cite{AT} (the last problem is whether the 
product of an $\mathbb R$-factorizable group and (a subgroup of) a $\sigma$-compact group is 
$\mathbb R$-factorizable). 

We use the notation $\mathbb R$ for the set of real numbers, $\mathbb N$ for the set of positive integers, 
and $\omega$ for the set of nonnegative integers. By $\oplus$ we denote the topological sum of spaces and   
by $|A|$, the cardinality of a set $A$. The definitions of the covering dimensions $\dim$ and $\dim_0$ can be 
found in \cite{Ch}. A topological space $X$ is \emph{zero-dimensional} if it has a base consisting of clopen sets 
and \emph{strongly zero-dimensional} if any finite cover of $X$ by cozero sets has a disjoint finite refinement 
(that is, $\dim_0(X)=0$). A subset $Y$ of a space $X$ is said to be \emph{$C$-embedded} in $X$ if any real-valued 
continuous function on $Y$ has a continuous extension to $X$, and $Y$ is \emph{$z$-embedded} in $X$ if every 
zero set of $Y$ is the trace on $Y$ of some zero set of $X$. 

The main result of this paper is the following theorem. 

\begin{theorem*}
There exist Boolean (and hence Abelian) Hausdorff topological groups $G_1$ and $G_2$ with the following properties:
\begin{enumerate}
\item[\textup{(i)}]
$G_1$ and $G_2$ are $\mathbb R$-factorizable and zero-dimensional;
\item[\textup{(ii)}] 
$G_1$ is submetrizable, and $G_1^n$ is Lindel\"of for any $n\in \mathbb N$;
\item[\textup{(iii)}]
$G_2$ is second countable;
\item[\textup{(iii)}]
$G_1\times G_2$ is not $\mathbb R$-factorizable.
\end{enumerate}
\end{theorem*}

Our construction of the groups $G_1$ and $G_2$ is based on 
Przymusi\'nski's notion of $n$-cardinality \cite{Prz1} and on his construction 
of Lindel\"of spaces $X$ and $Y$ such that $X\times Y$ is normal and $\dim X=\dim Y= 0$ but $\dim (X\times Y) >0$ 
\cite{Prz2}. Below we recall some details, following the exposition of the construction given in~\cite{Ch}.

\begin{definition}[\cite{Prz1}]
Let $X$ be a set, and let $n\in \mathbb N$. The \emph{$n$-cardinality} (with respect to $X$) of a set $A\subset 
X$, denoted by $|A|_n$, is the least cardinal $\kappa$ such that
$$
A\subset \bigcup_{i=1}^n (X^{i-1}\times Y\times X^{n-i})
$$
for some $Y \subset X$ with $|Y|=\kappa$ (here and in what follows it is assumed that $X^0\times Y=Y\times 
X^0=Y$). Clearly, $|A|_1= |A|$ and $|A|_n \le |A|$. If $|A|_n \le \omega$, then $A$ is 
said to be \emph{$n$-countable}; otherwise, $A$ is said to be \emph{$n$-uncountable}. 
\end{definition}

For $x\in X^n$ and $i\le n$, we denote the $i$th coordinate of $x$ by $x_i$ and the set of 
all coordinates of $x$ by $\tilde x$; in other words, $\tilde x= \{x_1, \dots, x_n\}$.

\begin{lemma}[{see \cite[Lemma 24.1]{Ch}}]
\label{lem1} 
Given a set $X$, a positive integer $n$, and an infinite cardinal $\kappa$, 
the following conditions on $A \subset X^n$ are equivalent:
\begin{enumerate}
\item[\textup{(a)}]
$|A|_n \ge \kappa$; 
\item[\textup{(b)}]
$A$ contains a subset $B$ of cardinality $\kappa$ such that $\tilde p\cap \tilde q =\varnothing$ whenever 
$p$ and $q$ are distinct points of $B$.
\end{enumerate}
\end{lemma}

\begin{definition}[{\cite[p.~186]{Ch}}]
Suppose given $n\in \mathbb N$, a set $X$, and a topology $\tau$ on $X^n$. A set $B 
\subset X$ is said to be \emph{weakly $n$-Bernstein with respect to} $\tau$ if $|A\cap B^n|_n=2^\omega$ 
for every $n$-uncountable $\tau$-closed set $A\subset X^n$. 
\end{definition}

If $\tau$ is a topology on a set $X$, then, abusing notation, we denote the product topology on~$X^n$ 
by $\tau^n$. The proof of the following lemma is very similar to that of Theorem~24.3 in~\cite{Ch}.

\begin{lemma}[{see \cite[Theorem 24.3]{Ch}}]
\label{lem2}
Let $(X,\tau)$ be a space with separable completely metrizable topology $\tau$, and let $\mu$ be a
topology on $X^2$ with the following properties: 
\begin{enumerate}
\item[\textup{(i)}]
$\mu\supset \tau^2$;
\item[\textup{(ii)}]
$X^2$ contains at most $2^\omega$ \,$2$-uncountable $\mu$-closed sets;
\item[\textup{(iii)}]
$|A|_2\ge 2^\omega$ for any $2$-uncountable $\mu$-closed set $A\subset X^2$.
\end{enumerate}
Then $X$ contains pairwise disjoint sets $B_1,B_2, \dots$ such that every $B_i$ is weakly $2$-Bernstein with 
respect to $\mu$ and weakly $n$-Bernstein with respect to $\tau^n$ for all~$n\in \mathbb N$.  
\end{lemma}  

\begin{proof}
Let us denote the family of all 
$2$-uncountable $\mu$-closed subsets of $X^2$ by $\mathscr  A_{\mu,2}$ and 
the family of all $n$-uncountable 
$\tau^n$-closed subsets of $X^n$, where $n\in \mathbb N$,  by $\mathscr  A_{\tau,n}$. 
Note that $\mathscr  A_{\mu,2}\supset \mathscr A_{\tau,2}$ (because $\tau^2\subset \mu$),  
$|\mathscr  A_{\tau,n}|\le 2^\omega$ for $n\in \mathbb N$ (because $(X^n,\tau^n)$ is second countable), 
and  $|\mathscr  A_{\mu,2}|\le 2^\omega$ (by assumption (ii) of the lemma).
We set 
$$
\mathscr  A=\mathscr A_{\mu,2}\cup\bigcup_{n\ne 2}\mathscr  A_{\tau,n}
$$ 
and index the elements of $\mathscr  A$ by ordinals less than $2^\omega$ as 
$\mathscr A=\{A_\alpha: \alpha< 2^\omega\}$ so that each 
element is assigned $2^\omega$ indices. Let $\alpha <2^\omega$. If 
$A_\alpha\in \mathscr  A_{\mu,2}$, then we set $n(\alpha)=2$; otherwise, 
we denote  by $n(\alpha)$ the unique $n\in \mathbb N$ ($n\ne 2$) for 
which $A_\alpha\in \mathscr  A_{\tau,n}$. Then,  for all 
$\alpha\in 2^\omega$ and $i\in \mathbb N$,
we recursively choose points $p(\alpha,i)\in A_\alpha$  so 
that $\tilde p(\alpha,i)\cap \tilde p(\beta,j)=\varnothing $ if $\alpha\ne \beta$ or $i\ne j$ 
in precisely the same way as in the proof 
of Theorem~24.3 of \cite{Ch}; the only difference is that, in the case $n(\gamma)=2$, we use our assumption (iii) 
and Lemma~\ref{lem1} to find a $B\subset A_\gamma$ such that $|B|=2^\omega$ 
and $\tilde p\cap \tilde q =\varnothing $ for any distinct $p,q\in B$. After that, following \cite[Theorem~24.3]{Ch}, 
we set 
$$
B_i=\bigcup \{\tilde p(\alpha,i):\alpha<2^\omega\}
$$ 
for each $i\in \mathbb N$. Clearly, $B_i\cap B_j=\varnothing$ if $i\ne j$. For each $n\ne 2$, any $n$-uncountable 
$\tau^n$-closed subset $A$ of $X^n$ equals $A_\alpha$ for $2^\omega$ indices $\alpha\in 2^\omega$, and we 
have $p(\alpha,i)\in A\cap B^{n(\alpha)}_i$ and $n(\alpha)=n$ for each of these $\alpha$ and all 
$i\in \mathbb N$. Since $\tilde p(\alpha, i)\cap \tilde p(\beta, i)=\varnothing$ for $\alpha\ne \beta$, 
it follows that $|A\cap B_i^n|_n\ge 2^\omega$ by Lemma~\ref{lem1}. 
Similarly, we have $|A\cap B_i^2|_2\ge 2^\omega$ for any $2$-uncountable 
$\mu$-closed (and hence for any $2$-uncountable $\tau^2$-closed) subset $A$ of $X^2$.
\end{proof}

Let $C$ be the Cantor set in $[0,1]\subset \mathbb R$. In \cite{Ch}
a topology $\mu$ on $C^2$ was defined which satisfies conditions (i)--(iii) of Lemma~\ref{lem2} for $C$ playing 
the role of $X$ and the  
usual (Euclidean) topology $\varepsilon$ on $C$ playing the role of $\tau$ 
(see the proofs of Lemmas~27.2 and 27.3 in \cite{Ch}). By 
Lemma~\ref{lem2} $C$ contains pairwise disjoint sets $S$, $S_1$, and $S_2$ which are weakly 
$2$-Bernstein with respect to $\mu$ and weakly $n$-Bernstein with respect to $\tau^n$ for all~$n\in \mathbb N$.  
In \cite[proof of Theorem~27.5]{Ch}, given any such sets, topologies 
$\tau_1$ and $\tau_2$ on $C$ were constructed which satisfied, in particular, the following conditions for $i=1,2$ 
(see~\cite[pp.~210, 211]{Ch}):
\begin{enumerate}
\item[(1)]
$\tau_i\supset \varepsilon$; 
\item[(2)]
any $\tau_i$-neighborhood of any point of $S_i$ is an $\varepsilon$-neighborhood of this point; 
\item[(3)]
$\tau_i$ has a base consisting of $\varepsilon$-closed sets;
\item[(4)]
$\dim(C, \tau_i)=\dim_0(C, \tau_i)=0$; 
\item[(5)] 
$(C,\tau_1)\times (C,\tau_2)$ is normal and $\dim((C,\tau_1)\times (C,\tau_2))=\dim_0((C,\tau_1)\times 
(C,\tau_2))=1$. 
\end{enumerate} 

We fix topologies $\tau_1$ and $\tau_2$ on $C$ with these properties and set 
$C_i=(C, \tau_i)$ for $i=1,2$. Note that it follows from (2) that the restriction of the topology 
$\tau_2$ to $S_2$ coincides with the topology induced on $S_2$ by $\varepsilon$. In what follows, by $S_2$ we 
mean the set $S_2$ endowed with this topology, i.e., treat $S_2$ as a subspace of $(C,\varepsilon)$; this is 
a separable metrizable space. In 
\cite[Example~27.8]{Ch} it was shown that 
\begin{enumerate}
\item[(6)]
$\dim (C_1 \times S_2) \ge \dim_0 (C_1 \times S_2)=1$.
\end{enumerate}

\begin{lemma}
\label{lem3}
The spaces $C_1^n$ are Lindel\"of for all $n\in \mathbb N$. 
\end{lemma}

\begin{proof}
We argue by induction on $n$.

Let $\gamma$ be a $\tau_1$-open cover of $C_1$. In view of (2), each point $s\in S_1$ has an $\varepsilon$-open 
neighborhood $U_s$ contained in an element $V_s$ of $\gamma$. Let $U=\bigcup_{s\in S_1} U_s$. Since $S_1$ is 
weakly $1$-Bernstein with respect to $\varepsilon$ and $C\setminus U$ is an $\varepsilon$-closed set disjoint 
from $S_1$, it follows that $C\setminus U$ is $1$-countable, that is, countable. For each $x\in C\setminus U$, 
choose an element $V_x$ of $\gamma$ containing $x$. Let $\{U_{s_k}:k\in \mathbb N\}$ be a countable subcover 
of the $\varepsilon$-open cover $\{U_s:s\in S_1\}$ of $S_1$. Then $\{V_{s_k}:k\in \mathbb N\}\cup 
\{V_x:x\in C\setminus U\}$ is a countable subcover of $\gamma$. 

Suppose that $n>1$ and $C_1^k$ is known to be Lindel\"of for every $k<n$. Let $\gamma$ be a 
$\tau_1^n$-open cover of $C_1^n$. In view of (2) each point $s\in S_1^n$ has 
an $\varepsilon^n$-open neighborhood $U_s$ contained in an element $V_s$ of $\gamma$. Let $U=\bigcup_{s\in S_1^n} 
U_s$. Since $S_1$ is 
weakly $n$-Bernstein with respect to $\varepsilon^n$ and $C_1^n\setminus U$ is an $\varepsilon^n$-closed set 
disjoint from $S_1^n$, it follows that $C_1^n\setminus U$ is $n$-countable, that is, there exists a countable 
set $Y\subset C_1$ such that 
$$
C_1^n\setminus U \subset \bigcup_{k=1}^n (C_1^{k-1}\times Y\times C_1^{n-k}).
$$
This means that $C_1^n\setminus U$ is contained in the countable union of spaces of the form 
$C_1^{k-1}\times \{x\}\times C_1^{n-k}$, where $k\le n$ and $x\in Y$, 
each of which is homeomorphic to $C_1^{n-1}$ and therefore Lindel\"of by the induction hypothesis. 
It remains to choose a countable subfamily of $\gamma$ 
covering $C_1^n\setminus U$ and a countable subfamily of $\{V_s:s\in S_1^n\}$ covering $U$, which exists because  
$\{V_s:s\in S_1^n\}$ has the $\varepsilon^n$-open refinement $\{U_s:s\in S_1^n\}$. 
\end{proof}

\begin{lemma}
\label{lem4}
Suppose that $G_1$, $G_2$, $M_1$, and $M_2$ are topological groups with the following 
properties:
\begin{enumerate}
\item[\textup{(i)}]
$M_1$ and $M_2$ are topological products of zero-dimensional second countable topological groups;
\item[\textup{(ii)}]
$G_1$ and $G_2$ are closed subgroups of $M_1$ and $M_2$, respectively; 
\item[\textup{(iii)}]
$C_1\times S_2$ is $C$-embedded in $G_1\times G_2$.
\end{enumerate}
Then the group $G_1\times G_2$ is not $\mathbb R$-factorizable. 
\end{lemma}

\begin{proof}
Any product of zero-dimensional 
second countable topological spaces is strongly zero-dimensional \cite{Morita}. Therefore, so is the product 
$M_1\times M_2$, and it contains $G_1\times G_2$ as a subgroup. As is known, any $\mathbb R$-factorizable subgroup 
of a topological group $G$ is $z$-embedded in $G$ \cite[Theorem~3.2]{HST}. It follows that if the group 
$G_1\times G_2$ were $\mathbb R$-factorizable, then this group, as well as its $C$-embedded 
subspace $C_1\times S_2$, would be $z$-embedded in $M_1\times M_2$. On the other hand, any $z$-embedded subspace 
of a strongly zero-dimensional space is strongly zero-dimensional \cite[Theorem~11.22]{Ch}, while 
$\dim_0(C_1\times S_2)>0$. Hence $C_1\times S_2$ is not $z$-embedded in $M_1\times M_2$ and 
$G_1\times G_2$ is not $\mathbb R$-factorizable.
\end{proof}

The product $C_1\times S_2$ is surely $C$-embedded in $G_1\times G_2$ if $C_1\times S_2$ is a retract of 
$G_1\times G_2$, which is the case if $C_1$ is a retract of $G_1$ and $S_2$ is a retract of $G_2$. 
Thus, we will look  for topological groups $G_1$ and $G_2$ containing $C_1$ and $S_2$ as retracts. 
These $G_1$ and $G_2$ will be 
the Boolean groups $B(C_1)$ and $B(S_2)$, respectively, with special topologies. 

A Boolean group is a group in which all elements are of order 2 (all such groups are Abelian), and 
the Boolean group $B(X)$ with basis $X$ is the set $[X]^{<\omega}$ of finite subsets of $X$ 
endowed with the operation $\vartriangle$ of symmetric difference. The zero element is the empty set. 
Each point $x\in X$ is identified with the 
singleton $\{x\}$. We use the notation $\+$ for the group operation of $B(X)$ and occasionally write 
$\vartriangle$ instead of $\+$. Thus, if $x\in X$, $F, G\in B(X)$, and $\mathbf A\subset B(X)$, then
\begin{gather*}
x\+F=\{x\}\+ F=\{x\}\vartriangle F, \qquad F\+G=F\vartriangle  G, \\
F\+\mathbf A = 
\{F\+ A: A\in  \mathbf A\}=\{F\vartriangle A: A\in  \mathbf A\}.
\end{gather*}
Let $X$ be a topological space. The subgroups of $B(X)$ of the form
$$
\mathbf{H}_\gamma=\{F\in B(X): \text{$|F\cap U|$ is  even for each $U\in \gamma$}\}, 
$$
where $\gamma$ ranges over all disjoint open covers of $X$, are normal (since $B(X)$ is Abelian), and 
the set of all these subgroups is obviously closed under the formation of finite intersections. Therefore, 
this set is a neighborhood base at zero of 
a group topology on $B(X)$ (see, e.g., \cite[Theorem~1.3.12]{AT}). If $X$ is zero-dimensional, 
then $B(X)$ with this topology contains $X$ as a subspace, because given any $\gamma$ and any $x\in X$, 
we obviously have $x\+\mathbf{H}_\gamma\cap X= U$, where $U$ is the element of $\gamma$ containing $X$.  
In what follows, we use the notation $B(X)$ for the abstract (that is, 
without topology) Boolean group with basis $X$ and $B^{\text{lin}}(X)$ 
for $B(X)$ with this topology. 

Recall that a topological space is said to be \emph{non-Archimedean} if it has a base such that, given  
any two of its elements, either they are disjoint or one of them contains the another (see~\cite{Nyikos}). 
In Theorem~3 (version~2) of \cite{3}, for a non-Archimedean space $X$, a retraction of the subspace 
$$
B_{\text{odd}}(X)=\{F\in B(X): \text{$|F|$ is odd}\}
$$ 
of $B^{\text{lin}}(X)$ onto $X$ was constructed (in \cite{3} the group $B^{\text{lin}}(X)$ was denoted by 
$B_z(X)$). In the particular case of the Cantor 
set $C$, the construction can be described as follows. 

Recall that $C$ can be represented as the subset of $[0,1]$ consisting of those numbers whose ternary expansions 
do not contain 1. This suggests the natural base $\mathscr B$ for the topology of $C$: 
$$
\mathscr B=\{U_{n_1\dots n_k}: k\in \mathbb N, n_i\in \{0,2\} \text{ for } i\le k\},
$$
where $U_{n_1\dots n_k}$ denotes the set of all numbers in $[0,1]$ whose ternary expansions begin with 
$0.n_1\dots n_k$. We also include the whole set $C$ in $\mathscr B$. Clearly, the elements of $\mathscr B$ 
form a tree with respect to reverse inclusion and every element of $\mathscr B$ is clopen. 

There are two natural orders on the set of subsets of $C$, the order by inclusion and the order induced by the 
usual order of $\mathbb R$. In what follows, when writing, say, ``$A < B$,'' ``$\min A$,'' or ``$A$ is on the 
left of $B$,'' we always mean the latter, unless otherwise is explicitly stated. Note that, given any two elements 
of $\mathscr B$, either one of them is contained in the other or one of them is on the left of the other. 

Let $F$ be any finite subset of $C$. We say that a set $A\subset C$ is \emph{$F$-void} if $A\cap  
F=\varnothing$, \emph{$F$-even} if $|A\cap F|$ is even and positive, and \emph{$A$-odd} if $|A\cap 
F|$ is odd. Clearly, each $F$-even element of $\mathscr B$ is contained in a maximal (by inclusion) $F$-even element 
of $\mathscr B$, and the union of these maximal $F$-even elements is equal to the union of all $F$-even elements 
of~$\mathscr B$. Moreover, this union is itself $F$-even, because $\mathscr B$ is a tree and therefore any two 
maximal $F$-even elements either coincide or are disjoint. 

For $F\in B_\text{odd}(C)$, we set   
$$
r(F)=\min(F\setminus \textstyle{\bigcup} \{B\in \mathscr B: \text{$B$ is $F$-even}\}), \eqno(*)
$$
or, equivalently,
\begin{multline*}
r(F)= \min(F\setminus \textstyle{\bigcup} \{B\in \mathscr B:\\ 
\text{$B$ is an inclusion-maximal $F$-even  element of 
$\mathscr B$}\}).
\end{multline*}

\begin{lemma}
\label{lem5}
There exists a second countable zero-dimensional group topology $\tau$ on $B(C)$ such that 
it induces the Euclidean topology $\varepsilon$ on $C$ and the map $r\colon B_{\text{\rm odd}}(C)\to C$ 
defined by $(*)$
is continuous with respect to $\tau$ restricted to $B_{\text{\rm odd}}(C)$.
\end{lemma}

\begin{proof}
Recall that, given a disjoint open cover $\gamma$ of $C$, 
$$
\mathbf{H}_\gamma= \{F\in B(C): \text{$|F\cap U|$ is  even for each $U\in \gamma$}\}.
$$
The family 
$$
\mathscr H=\{\mathbf{H}_\gamma :\text{$\gamma$ is a disjoint cover of $C$ by elements of $\mathscr B$}\}
$$
of subgroups is a neighborhood base at zero for a group topology $\tau$ of $B(C)$. This family is countable, 
because all open disjoint covers of $C$ are finite (since $C$ is compact) 
and $\mathscr B$ is countable. It is easy to check that $C$ is contained in $(B(C), \tau)$ as a 
subspace. Indeed, take  any point  $x\in C$ and any neighborhood $V_x\in \mathscr B$ of $x$. Let 
$\gamma$ be a disjoint cover of $C$ consisting of $V_x$ and some other elements of $\mathscr B$. 
If $F\in \mathbf{H}_\gamma$ and $x \+ F=\{x\}\vartriangle F\in C$, then either $F=\varnothing$ or $F=\{x,y\}$. In the 
latter case, $x\+F=y$ and by the definition of $\mathbf{H}_\gamma$ the point $y$ must belong to the same element of 
$\gamma$ as $x$, that is, to $V_x$. Thus, $(x\+\mathbf{H}_\gamma)\cap C\subset V_x$. This shows that the topology 
induced by $\tau$ on $C$ is not coarser than the topology $\varepsilon$ of $C$. On the other hand, it cannot be 
finer, because $\tau$ is coarser than the topology of $B^{\text{lin}}(C)$. 
Obviously, $(B(C),\tau)$ is $T_0$ and hence Tychonoff.  

Note that all elements in any $\mathbf{H}_\gamma$ are of even cardinality.  Therefore, for every $F\in B_\text{odd}(C)$, we 
have $F\+ \mathbf{H}_\gamma=\{F\vartriangle H: H\in \mathbf{H}_\gamma\}\subset B_\text{odd}(C)$. 

Let us show that the map $r$ is continuous with respect to the topology $\tau|_{B_\text{odd}(C)}$. Suppose that 
$x=r(F)$ for $F\in B_\text{odd}(X)$. By construction $x\in F$. Take any neighborhood $U$ of $x$. Let $V_1,\dots, V_m$ be 
all inclusion-maximal $F$-even elements of $\mathscr B$; their number is finite because they are pairwise 
disjoint (since $\mathscr B$ is a tree) and each of them intersects the finite set $F$. None of 
these elements contains $x$ (because $x=r(F)$), and all of  them are clopen. Choose a neighborhood $V_x\in 
\mathscr B$ of $x$ satisfying the conditions $V_x\subset U$, $V_x\cap F=\{x\}$, and $V_x\cap V_i=\varnothing$ for $i\le 
m$. Consider the cover of $C$ consisting of the sets $V_x$ and $V_1,\dots, V_m$ and of all elements of $\mathscr B$ 
disjoint from them. This cover has a disjoint subcover $\gamma$, because any two of its elements are either 
disjoint or contained in one another (recall that $\mathscr B$ is a tree). Clearly, $\gamma$ is finite. 
We claim that $r(F\+\mathbf{H}_\gamma)\subset V_x$. 

Indeed, take an $H\in \mathbf{H}_\gamma$. We must show that $r(F\vartriangle H)\in V_x$. 
Note that an element $V$ of $\gamma$ is ($F\vartriangle H$)-odd if and only if it is $F$-odd, 
because each element of $\gamma$ is either $H$-even or $H$-void and a point of $F$ can be 
cancelled in $F\vartriangle H$ only by some point of $H$. 

Let $V$ be the leftmost (with respect to the natural order $<$ on $C$) $F$-odd (=~($F\vartriangle H$)-odd) 
element of $\gamma$. Note that $V\cap F$ is disjoint from all 
inclusion-maximal $F$-even elements of $\mathscr B$, because all such elements are included in  
$\gamma$ and $V$ is not among them. By the definition of the map $r$ we have $x=r(F)\le \min (V\cap F)$. Since 
$x\in V_x$, it follows that $V_x$ either coincides with $V$ or is on the left of $V$, and since $V_x$ is $F$-odd, 
it follows that $V=V_x$. 

The point $r(F\vartriangle H)$ cannot belong to an ($F\vartriangle H$)-even or ($F\vartriangle H$)-void element of $\gamma$, because 
$\gamma\subset \mathscr B$ and 
\begin{align*}
r(F\vartriangle H)&\in (F\vartriangle H)\setminus \bigcup \{B\in \mathscr B: \text{$B$ is $(F\vartriangle H)$-even}\}\\
&=(F\vartriangle H)\setminus \bigcup \{B\in \mathscr B: \text{$B$ is $(F\vartriangle H)$-even or $(F\vartriangle H)$-void}\}).
\end{align*}
Therefore, the element $V$ of $\gamma$ containing $r(F\vartriangle H)$ is ($F\vartriangle H$)-odd and hence either 
coincides with $V_x$ or is on the right of $V_x$. Since $r(F\vartriangle H)$ is the least element of $F\vartriangle H$ 
not belonging to $\bigcup \{B\in \mathscr B:$ $B$ is $(F\vartriangle H)$-even$\}$ and $r(F\vartriangle H)\in V$, it follows that 
there exists a family $\mathscr B'$ of $(F\vartriangle H)$-even elements of $\mathscr B$ such that 
$$
\bigcup\mathscr B'\supset\{y\in F\vartriangle H: y<V\}.
$$ 

Suppose that $V\ne V_x$. Note that $V_x$ is $(F\vartriangle H)$-odd, because it contains precisely one point of $F$ and 
an even number (or none) of the points of $H$. Hence $(F\vartriangle H)\cap V_x\ne \varnothing$. Let  
$W_1, \dots, W_k$ be all inclusion-maximal elements of $\mathscr B'$ intersecting $(F\vartriangle H)\cap V_x$. Each 
$W_i$, being an element of $\mathscr B$, either contains $V_x$ or is contained in $V_x$, 
because $\mathscr B$ is a tree, $V_x\in \mathscr B$, and $\mathscr B'\subset \mathscr B$. By maximality 
the sets $W_1, \dots,  W_k$ are pairwise disjoint. Therefore, if $k\ge 2$, then all of them are contained in 
$V_x$ and $V_x\cap (F\vartriangle H)= \bigcup_{i\le k} W_i\cap (F\vartriangle H)$. This is impossible,  because $|V_x\cap (F\vartriangle 
H)|$ is odd and  all $|W_i\cap (F\vartriangle H)|$ are even. Thus, some element $W$  of $\mathscr B'$ contains $V_x\ni x$. 
Moreover, this $W$ is a union of some elements of $\gamma$, since $W\in \mathscr B$, $\gamma\subset \mathscr B$, 
and $\gamma$ covers $C$. This means that $|W\cap H|$ is even and hence so is 
$|W\cap F|$, because $W$ is $(F\vartriangle H)$-even. However, $x$  equals $r(F)$ and hence does not belong to any 
$F$-even or $F$-void element of $\mathscr B$. This contradiction proves that $V=V_x$, i.e., $r(F\vartriangle H)\in V_x$. 

Thus, $r(F\+\mathbf{H}_\gamma)\subset V_x$. We have shown that, for any $F\in B_\text{odd}(C)$ and any neighborhood $U$ of 
$x=r(F)$ in $C$, there exists an $\mathbf{H}_\gamma\in \mathscr H$ such that the image of the open neighborhood $F\+ 
\mathbf{H}_\gamma$ of $F$ in $(B(C),\tau)$ under $r$ is contained in $U$. This means that $r$ is continuous with 
respect to the topology $\tau|_{B_\text{odd}(C)}$. 

It remains to note that the group $(B(C), \tau)$ is zero-dimensional and metrizable, because 
the neighborhood base $\mathscr H$ at zero is countable and consists of open (and hence closed) subgroups, 
and it is separable, 
because 
$$
B(C)=\bigcup_{n\in \omega} B_n(C),\quad \text{where} \quad 
B_n(C) = \{F\in B(C): |F|\le n\} 
$$
and each $B_{n}(C)$ is the image of the separable space $(C\oplus \{\varnothing\})^n$ under the 
addition map $i_n\colon (x_1,\dots, x_n)\mapsto x_1\+ \dots\+x_n$, which is 
continuous with respect to any group topology on $B(C)$ inducing $\epsilon$ on~$C$.  
\end{proof}

Let $B^\tau(C)$ denote the group $B(C)$ with the topology $\tau$ defined in Lemma~\ref{lem5}. 

\begin{lemma}
\label{lem6}
The Cantor space $C$ is a retract of $B^\tau(C)$. Moreover, for any $x_0\in C$, the map 
$$
\hat r=\begin{cases}
r(F) &\text{if $F\in B_\text{\rm odd}(C)$},\\
x_0&\text{otherwise}, 
\end{cases}
$$
where $r$ is defined by $(*)$, is a retraction. 
\end{lemma}

\begin{proof}
The map $\hat r$ is continuous, because $B_\text{odd}(C)$ is clopen in $B^\tau(C)$, being a coset of the 
open subgroup 
$$
B_\text{even}(C)=\mathbf{H}_{\{C\}}=\{F\in B(C): \text{$|F\cap C|=|F|$ is  even}\}
$$
of $B^\tau(C)$. Clearly, for every $x\in C$, we have 
$$
\hat r(x)=r(x)=\min(\{x\}\setminus \textstyle{\bigcup} \{B\in \mathscr B: \text{$B$ is $\{x\}$-even}\})=x, 
$$
because there are no $\{x\}$-even sets. Thus, $\hat r$ is a retraction.
\end{proof}

Now we can prove the main theorem.

\begin{proof}[Proof of the main theorem] 
We take the group $B^{\text{lin}}(C_1)$ as $G_1$ and 
the subgroup of $B^\tau(C)$ generated by $S_2$ as~$G_2$. 
According to \cite[Theorem~7 (version 2)]{3}, $C_1$ is a retract of $B^{\text{lin}}(C_1)$. 
Take any $x_0\in S_2$. Restricting the retraction $\hat r$ defined in Lemma~\ref{lem6} for this $x_0$ to $G_2$, 
we obtain a retraction 
of $G_2$ onto $S_2$. Indeed, according to $(*)$, we have $r(F)\in F$ for any $F\in B(C)$. Therefore, 
$\hat r(F)\in F\cup \{x_0\}\subset S_2$ for any $F\in G_2$, whence $\hat r(G_2)= S_2$.

By Lemma~\ref{lem5} the group $B^\tau(C)$ is second countable and zero-dimensional; hence so is its subgroup 
$G_2$. The topology of $B^{\text{lin}}(C_1)$ is finer than $\tau$, which implies the submetrizability of $G_1$. 
The same argument as at the end of the proof of Lemma~\ref{lem5} shows that $G_1^n$ is Lindel\"of for any 
$n\in \mathbb N$. In more detail, 
$$
B^{\text{lin}}(C_1)=\bigcup_{n\in \omega} B_n(C),\quad \text{where} \quad 
B_n(C_1) = \{F\in B(C): |F|\le n\},
$$
and each $B_{n}(C_1)$ is the image of $(C_1\oplus \{\varnothing\})^n$ under the continuous 
addition map $i_n\colon (x_1,\dots, x_n)\mapsto x_1\+ \dots\+x_n$. Hence $G_1=B^{\text{lin}}(C_1)$ is a continuous image 
of the sum $C_\infty=\bigoplus_{n\in \omega}(C_1\oplus \{\varnothing\})^n$ and 
$G_1^n$ is a continuous image of $C_\infty^n$ for every $n\in \mathbb N$. By Lemma~\ref{lem3} all spaces 
$C_1^n$ are Lindel\"of; therefore, so are $(C_1\oplus \{\varnothing\})^n$ and $C_\infty^n$. It follows that 
all $G_1^n$ are Lindel\"of. 

Note that both groups $G_1$ and $G_2$ are $\mathbb R$-factorizable, being Lindel\"of~\cite{Tk1}. Let us show that 
$G_1\times G_2$ is not. To this end, we first embed $G_1$ 
in a product of zero-dimensional second countable groups and then apply Lemma~\ref{lem4}.

Let $\Gamma$ denote the set of all disjoint open covers of $C_1$.  We fix a countable discrete space 
$D= \{d_n: n\in \mathbb N\}$ and denote by $B^{\text{d}}(D)$ the Boolean group with basis $D$ endowed 
with the discrete topology. Note that any cover $\gamma\in \Gamma$ is countable, because $C_1$ is Lindel\"of. 
Let $\gamma=\{U_n: n\in \mathbb N\}$ be such a cover. Consider the map $f_\gamma\colon C_1\to D$ defined by 
$f_\gamma(U_n)= \{d_n\}$ for $n\in \mathbb N$. Let  
$\hat f_\gamma\colon G_1\to B^{\text{d}}(D)$ be the homomorphism extending $f$ to $G_1$; it is defined by 
$\hat f_\gamma(x_1\+\dots\+ x_n) = f_\gamma(x_1)\+ \dots \+ f_\gamma(x_1)$ for $x_1,\dots,x_n\in C_1$. 
The preimage $\hat f_\gamma^{-1}$ of the zero element $\varnothing$ of $B^{\text{d}}(D)$ is precisely 
$\mathbf{H}_\gamma = \{F\in B^{\text{lin}}(C_1): \text{$|F\cap U|$ is  even for each $U\in \gamma$}\}$; therefore, 
$\hat f_\gamma$ is continuous. Since the subgroups $\mathbf{H}_\gamma$, $\gamma\in \Gamma$, form a base of neighborhoods 
of zero for the topology of $B^{\text{lin}}(X)$, it follows that the homomorphims $\hat f_\gamma$, 
$\gamma\in \Gamma$, separate points from closed sets and therefore the diagonal 
$$
\mathop{\large\Delta}\limits_{\gamma\in \Gamma}\hat f_\gamma\colon B^{\text{lin}}(C_1)\to 
B^{\text{d}}(D)^{|\Gamma|}
$$
is a homeomorphic embedding; clearly, this is a homomorphism. Thus, $B^{\text{lin}}(C_1)$ is topologically 
isomorphic to a subgroup of the product $B^{\text{d}}(D)^{|\Gamma|}$ and $B^{\text{d}}(D)$ is 
discrete and countable (and hence zero-dimensional and second countable), which is what we need. 

Applying Lemma~\ref{lem4} to the groups $G_1$, $G_2$, $M_1=B^{\text{d}}(D)^{|\Gamma|}$, and $M_2=G_2$, we see 
that $G_1\times G_2$ is not $\mathbb R$-factorizable.
\end{proof}

\medskip

The author is most grateful to Evgenii Reznichenko for very useful discussions.

\end{document}